\documentclass[11pt,reqno]{article}
\usepackage{euscript,amstext,amsthm,amssymb}
\usepackage{amsmath}

\theoremstyle{definition}

\theoremstyle{remark}

\numberwithin{equation}{section}

\newcommand{\be}{\begin{equation}}
      \newcommand{\ee}{\end{equation}}
      \newcommand{\ba}{\begin{eqnarray}}
       \newcommand{\ea}{\end{eqnarray}}
\newcommand{\ban}{\begin{eqnarray*}}
       \newcommand{\ean}{\end{eqnarray*}}

\newcommand{\pt}{\partial}

\newcommand{\ra}{\rightarrow}

 \newcommand{\Pf}{\noindent {\em Proof.} }

\newcommand{\bM}{\partial M}

\newcommand{\reta}{\overline{\eta}}

\newcommand{\Det}{\mbox{\rm Det}}

\newtheorem{theo}{Theorem}[section]

\begin{document}
\newtheorem{lem}[theo]{Lemma}
\newtheorem{prop}[theo]{Proposition}
\newtheorem{coro}[theo]{Corollary}
\newtheorem{rk}[theo]{Remark}

\title{Eta Invariant and Holonomy: the Even Dimensional Case}
\author{Xianzhe Dai\thanks{Math Dept, UCSB, Santa Barbara, CA 93106, USA \tt{Email:
dai@math.ucsb.edu}. Partially supported by NSF and NSFC.} \and
Weiping Zhang\thanks{Chern Institute of Mathematics and LPMC,
Nankai University, Tianjin 300071, P. R. China. \tt{Email:
weiping@nankai.edu.cn}. Partially supported by MOE and NNSFC} }

\maketitle
\begin{abstract}
In previous work, we introduced eta invariants for even
dimensional manifolds. It plays the same role as the eta invariant
of Atiyah-Patodi-Singer, which is for odd dimensional manifolds.
It is associated to  $K^1$ representatives on  even dimensional
manifolds and is closely related
to the so called WZW theory in physics. In fact, it is an intrinsic interpretation
of the Wess-Zumino term without passing to the bounding
$3$-manifold.  Spectrally the eta invariant is defined on a finite cylinder, rather than on the
manifold itself. Thus it is an interesting question to find an
intrinsic spectral interpretation of this new invariant. We
address this issue here using adiabatic limit technique. The
general formulation relates the (mod $\mathbb Z$ reduction of) eta
invariant for even dimensional manifolds with the holonomy of the determinant
line bundle of a natural family of Dirac type operators. In this
sense our result might be thought of as an even dimensional
analogue of Witten's holonomy theorem proved by Bismut-Freed
 and Cheeger independently.
\end{abstract}

\section{Introduction}


The $\eta$-invariant is introduced by Atiyah-Patodi-Singer in their
seminal series of papers \cite{aps1,aps2,aps3} as the correction
term from the boundary for the index formula on a manifold with
boundary. It is a spectral invariant associated to the natural
geometric operator on the (boundary) manifold and it vanishes for even
dimensional manifolds (in this case the corresponding manifold with
boundary will have odd dimension). In our previous work \cite{dz},
we introduced an invariant of eta type for even dimensional
manifolds. It plays the same role as the eta invariant of
Atiyah-Patodi-Singer.

 Any elliptic differential operator on an odd dimensional closed manifold will
have index zero. In this case, the appropriate index to consider is
that of Toeplitz operators. This also fits perfectly with the
interpretation of the index of Dirac operator on even dimensional
manifolds as a pairing between the even $K$-group and $K$-homology.
Thus in the odd dimensional case one considers the odd $K$-group and
odd $K$-homology. For a closed manifold $M$, an element of
$K^{-1}(M)$ can be represented by a differentiable map from $M$ into
the unitary group
\be \label{kor} g: \ M \longrightarrow { U}(N), \ee where $N$ is a positive integer. As we mentioned the
appropriate index pairing between the odd $K$-group and $K$-homology
is given by that of the Toeplitz operator, defined as follows.

Consider $L^2(S(TM)\otimes E)$,  the space of $L^2$ spinor
fields\footnote{In this paper,  for simplicity,  we will generally
assume that our manifolds are spin, although our discussion extends
trivially to the case of Dirac type opreatros.} twisted by an
auxilliary vector bundle $E$. It decomposes into an orthogonal
direct sum
$$L^2(S(TM)\otimes E)=\bigoplus_{\lambda\in {\rm Spec}(D^E)} E_\lambda,$$
according to the eigenvalues  $\lambda$ of the Dirac operator $D^E$. The ``Hardy space" will be
$$L^2_{\geq 0}(S(TM)\otimes E) =\bigoplus_{\lambda\geq 0}E_\lambda. $$
The corresponding orthogonal projection from $L^2(S(TM)\otimes E)$
to $L^2_{\geq 0}(S(TM)\otimes E)$ will be denoted by $P^E_{\geq 0}$.

The Toeplitz operator $T^E_g$ is then defined as \be \label{toe}
T^E_g=P^E_{\geq 0}gP^E_{\geq 0}:L^2_{\geq 0}\left(S(TM)\otimes E
\otimes {\bf C}^N\right) \longrightarrow L^2_{\geq
0}\left(S(TM)\otimes E \otimes {\bf C}^N\right).\ee This is a
Fredholm operator whose index is given by \be \label{tifcm} {\rm
ind}\, T^E_g=-\left\langle \widehat{A}(TM){\rm ch}(E){\rm
ch}(g),[M]\right\rangle ,\ee where ${\rm ch}(g)$ is the odd Chern
character associated to $g$ \cite{bd}. It is represented by the
differential form (cf. [Z1, Chap. 1])
$$
{\rm ch}(g)= \sum_{n=0}^{\dim M-1\over 2} {n!\over (2n+1)!}{\rm
Tr}\left[\left(g^{-1}dg\right)^{2n+1}\right].
$$

In \cite{dz} we established  an index theorem which generalizes
(\ref{tifcm}) to the case where $M$ is an odd dimensional spin
manifold with boundary $\partial M$. The definition of the
Toeplitz operator now uses Atiyah-Patodi-Singer boundary
conditions on $\partial M$. The self adjoint Atiyah-Patodi-Singer
boundary conditions depend on choices of Lagrangian subspaces $L
\subset \ker D^E_{\partial M}$. We will denote the corresponding
boundary condition by $P_{\partial M}(L)$. The resulting Toeplitz
operator will then be denoted by $T^E_g(L)$.

We recall the main result in \cite{dz} as follows.

\begin{theo} \label{dz} The Toeplitz operator $T^E_g(L)$ is Fredholm with index given
by \begin{eqnarray}  \label{titfmwb}
 {\rm ind}\, T^E_g(L)& = & -\left({1\over 2\pi\sqrt{-1}}\right)^{(\dim M+1)/2}\int_M
\widehat{A}\left(R^{TM}\right){\rm
Tr}\left[\exp\left(-R^E\right)\right]{\rm ch}(g) \nonumber \\
& & - \ \overline{\eta}(\partial M, E, g) + \tau_\mu
\left(gP_{\partial M}(L) g^{-1},P_{\partial M}(L) ,  {\mathcal
P}_M\right) .
\end{eqnarray}
\end{theo}

 Here $\overline{\eta}(\partial M, E, g)$ denotes the
invariant of $\eta$-type for even dimensional manifold $\pt M$ and
the $K^1$ representative $g$. The third term is an interesting new
{\em integer} term, a triple Maslov index introduced in [KL]. See
\cite{dz} for details.

\noindent{\bf Remark}. Our index formula is closely related
to the so called WZW theory in physics \cite{w0}. When $\bM=S^2$ or
a compact Riemann surface and $E$ is trivial, the local term in
(\ref{dz}) is precisely the Wess-Zumino term, which allows an integer
ambiguity, in the WZW theory. Thus, our eta invariant
$\overline{\eta}(\partial M, g)$ gives an intrinsic interpretation
of the Wess-Zumino term without passing to the bounding
$3$-manifold. In fact, for $\bM=S^2$, it can be further reduced to
a local term on $S^2$ by using Bott's periodicity, see \cite[Remark
5.9]{dz}.

The eta invariant $\overline{\eta}(\partial M, E, g)$ is defined
on a finite cylinder $[0, 1]\times \pt M$, rather than on $\pt M$
itself. Thus it is an interesting question to find an intrinsic
spectral interpretation of this new invariant. In this paper we
answer this question by using the adiabatic limit technique.
First, under invertibility assumptions, we give an explicit
formula for our eta invariant in terms of a natural family of
Dirac type operators on the manifold. This family arises from the
original Dirac operator by a perturbation involving the $K^1$
representative. The general formulation relates (the mod $\mathbb
Z$ reduction of) the eta invariant for even dimensional manifolds
with the holonomy of the determinant line bundle of this natural
family of Dirac type operators. The work of \cite{d} on the
adiabatic limit of eta invariants for manifolds with boundary and
that of \cite{df} on Witten's Holonomy Theorem play an important
role here.

This paper is organized as follows. In Section 2, we review the
definition of the eta invariant for an even dimensional closed
manifold introduced in \cite{dz}. In Section 3, we give an intrinsic
spectral interpretation of the eta invariant under certain
invertibility assumption. Section 4 deals with the general case.  And
we end with a conjecture and a few remarks in the last section.

Some of the results in this paper have been described in \cite{d3}.

\section{ An invariant of $\eta$ type for even dimensional
manifolds}

For an even dimensional closed manifold $X$ (which may or may not be
the boundary of an odd dimensional manifold) and a $K^1$
representative $g: X \ra U(N)$, the eta invariant will be defined in
terms of an eta invariant on the cylinder $[0, 1]\times X$ with
appropriate APS boundary conditions.

In general, for a compact manifold $M$ with boundary $\pt M$ with
the product structure near the boundary, the Dirac operator $D^E$
twisted by an hermitian vector bundle $E\otimes{\bf C}^N$
decomposes near the boundary as \be \label{ps} D^E =
c\left({\partial \over
\partial x}\right) \left( {\partial \over
\partial x} + D^E_{\partial M} \right). \ee
The APS projection $P_{\partial M}$ is an elliptic global boundary
condition for $D^E$. However, for self adjoint boundary
conditions, we need to modify it by a Lagrangian subspace of $\ker
D^E_{\partial M}$, namely, a subspace $L$ of $\ker D^E_{\partial M}$
such that $c({\partial \over
\partial x})L=L^\perp \cap (\ker D^E_{\partial M})$. Since $\partial
M$ bounds $M$, by the cobordism invariance of the index, such
Lagrangian subspaces always exist.

The modified APS projection is then obtained by reducing the kernel
part of the projection to the projection onto the Lagrangian
subspace. More precisely, denote
$$L^2_{+}((S(TM) \otimes E\otimes{\bf
C}^N)|_{\partial M}) =\bigoplus_{\lambda>
0}E_\lambda(D^{E\otimes{\bf C}^N}_{\partial M}), $$ where $\lambda$
runs over the positive eigenvalues of $D^{E\otimes{\bf
C}^N}_{\partial M}$. Denote by $P_{\partial M}$ the orthogonal
projection from $L^2((S(TM) \otimes E\otimes{\bf C}^N)|_{\partial
M})$ to $L^2_{+}((S(TM)\otimes E\otimes{\bf C}^N)|_{\partial M})$.
 Let $P_{\partial M}(L)$ denote the orthogonal
projection operator from $L^2((S(TM) \otimes E\otimes{\bf
C}^N)|_{\partial M})$ to $L^2_{+}((S(TM)\otimes E\otimes{\bf
C}^N)|_{\partial M})\oplus L$: \be \label{mapsp} P_{\partial M}(L)
=P_{\partial M}+P_L, \ee where $P_L$ denotes the orthogonal
projection  from $L^2((S(TM)\otimes E\otimes{\bf C}^N)|_{\partial
M})$ to $L$.

The pair $(D^E, P^E_{\partial M}(L) )$ forms a self-adjoint elliptic
boundary problem, and $P_{\partial M}(L) $ is called an
Atiyah-Patodi-Singer boundary condition associated to $L$. We will
denote the corresponding elliptic self-adjoint operator by
$D^E_{P_{\partial M}(L) }$.

In \cite{dz}, we originally intend to consider the conjugated
elliptic boundary value problem $D^{E}_{gP_{\partial M}(L)
g^{-1}}$ (cf. \cite{z2}). However, the analysis turns out to be
surprisingly subtle and difficult. To circumvent this difficulty,
a perturbation of the original problem was constructed.

Let $\psi=\psi(x)$ be a cut off function which is identically $1$ in
the $\epsilon$-tubular neighborhood of $\partial M$ ($\epsilon
>0$ sufficiently small) and vanishes outside the $2\epsilon$-tubular
neighborhood of $\partial M$. Consider the Dirac type operator
$$
D^{\psi}=(1-\psi)D^E + \psi gD^E g^{-1}.
$$

The motivation for considering this perturbation is that, near the boundary, the
operator $D^{\psi}$ is actually given by the conjugation of $D^E$,
and therefore, the elliptic boundary problem $(D^{\psi},
gP_{\partial M}(L) g^{-1})$ is now the conjugation of the APS
boundary problem $(D^E, P_{\partial M}(L))$, i.e., this is now
effectively standard APS situation and we have a self adjoint
boundary value problem $(D^{\psi}, gP_{\partial M}(L) g^{-1})$
together with its associated self adjoint elliptic operator
$D^{\psi}_{gP_{\partial M}(L) g^{-1}}$.

The same thing can be said about the conjugation of $D^{\psi}$: \be
\label{pdo} D^{\psi, g}=g^{-1} D^{\psi} g =D^E + (1-\psi)
g^{-1}[D^E, g] . \ee
We will in fact use $D^{\psi, g}$.

We are now ready to construct the eta invariant for even
dimensional manifolds. Given an even dimensional closed spin
manifold $X$, we consider the cylinder $[0, 1]\times X$ with the
product metric. Let $g: \ X \ra U(N)$ be a map from $X$ into the
unitary group which extends trivially to the cylinder. Similarly,
$E \ra X$ is an Hermitian vector bundle which is also extended
trivially to the cylinder. We assume that ${\rm ind}\, D^E_+=0$ on
$X$ which guarantees the existence of the Lagrangian subspaces
$L$.

Consider the analog of $D^{\psi,g}$ as defined in (\ref{pdo}), but
now on the cylinder $[0, 1]\times X$  and denote it by
$D^{\psi,g}_{[0, 1]}$. Here  $\psi=\psi(x)$ is a cut off function on $[0, 1]$ which is identically $1$ for  $0\leq x \leq \epsilon$ ($\epsilon
>0$ sufficiently small) and vanishes when $1-2\epsilon\leq x \leq 1$.  We equip it with the boundary condition
$P_{X}(L) $ on one of the boundary components $\{0\}\times X$ and
the boundary condition ${\rm Id}-g^{-1}P_{X}(L) g $ on the other
boundary component $\{1\}\times X$ (Note that the Lagrangian
subspace $L$ exists by our assumption of vanishing index). Then
$(D^{\psi,g}_{[0, 1]}, P_X(L) , {\rm Id}-g^{-1}P_{X}(L) g)$ forms a
self-adjoint elliptic boundary problem. For simplicity, we will
still denote the corresponding elliptic self-adjoint operator by
$D^{\psi, g}_{[0, 1]}$.

Let $\eta(D^{\psi, g}_{[0, 1]},s)$ be the $\eta$-function of
$D^{\psi, g}_{[0, 1]}$ which, when ${\rm Re}(s)>>0$, is defined by
\be \label{ef} \eta(D^{\psi, g}_{[0, 1]},s) =\sum_{\lambda \neq
0}{{\rm sgn}(\lambda)\over |\lambda|^s}, \ee
 where $\lambda$ runs
through the nonzero eigenvalues of $D^{\psi, g}_{[0, 1]}$.

By \cite{mu,df}, one knows that the $\eta$-function $\eta(D^{\psi,
g}_{[0, 1]},s)$ admits a meromorphic extension to ${\bf C}$ with
$s=0$ a regular point (and it has only simple poles). One then
defines, as in \cite{aps1}, the $\eta$-invariant of
$D^{\psi,g}_{[0,1]}$ by $\eta(D^{\psi,g}_{[0, 1]})= \eta
(D^{\psi,g}_{[0, 1]},0 )$,  and the reduced $\eta$-invariant by
\be \label{rei} \overline{\eta}\left(D^{\psi,g}_{[0,
1]}\right)={\dim \ker D^{\psi,g}_{[0, 1]} +
\eta\left(D^{\psi,g}_{[0, 1]}\right) \over 2}. \ee

\noindent {\bf Definition 2.1.} {We define an invariant of $\eta$
type for the Hermitian  vector bundle $E$ on the even dimensional
manifold $X$ (with vanishing index) and the $K^1$ representative
$g$ by \be \label{eifedm} \overline{\eta}(X,E, g)=
 \overline{\eta} \left(D^{\psi, g}_{[0, 1]}\right) - {\rm sf} \left\{D^{\psi,g}_{[0, 1]}(s);
 0 \leq s \leq 1 \right\},
\ee where $D^{\psi,g}_{[0, 1]}(s)$ is a path connecting $g^{-1} D^E
g$ with $D^{\psi,g}_{[0, 1]}$ defined by
$$
D^{\psi,g}(s)= D^E + (1-s\psi) g^{-1}[D^E, g]
$$
on $[0, 1]\times X$, with the boundary condition $P_{X}(L) $ on
$\{0\}\times X$ and the boundary condition ${\rm
Id}-g^{-1}P_{X}(L) g$ at $\{1\}\times X$.}
\newline

 It was shown in \cite{dz} that $\overline{\eta}(X, E,
g)$ does not depend on the cut off function $\psi$.

\section{An intrinsic spectral interpretation, the invertible case}

The usefulness of the eta invariant of Atiyah-Patodi-Singer comes,
at least partially, from the spectral nature of the invariant,
i.e. that it is defined via the spectral data of the Dirac
operator on the (odd dimensional) manifold. Our eta invariant for
even dimensional manifold is defined via the eta invariant on the
corresponding odd dimensional cylinder by imposing APS boundary
conditions. Thus, it will be desirable to have a direct spectral
interpretation in terms of the spectral data of the original
manifold (and the $K^1$ representative). In this section we give
such an interpretation under certain invertibility assumption.
This invertibility condition will be removed in the next section.

The crucial point here is the following observation. As in the
previous section, we can also consider the invariant
$\overline{\eta}(D^{\psi,g}_{[0, a]})$, similarly constructed on a
cylinder $[0, a] \times X$ of radial size $a>0$.

\begin{lem} Assuming that $ \ker [D_X + s\ c(g^{-1} dg)] =0,\ \ \forall\  0
\leq s \leq 1$, then $\overline{\eta}(D^{\psi,g}_{[0, a]})$, and
hence $\overline{\eta}(X,E, g)$, is independent of $a$. Without
the invertibility assumption, the mod $\mathbb Z$ reduction of
$\overline{\eta}(X,E, g)$ is independent of $a$.
\end{lem}

This can be seen by a rescaling argument (cf. [M\"u, Proposition
2.16], see also \cite[Theorem 3.2]{d2}).

On the other hand, as we mentioned before,
\begin{lem}
$\overline{\eta}(X,E, g)$, is independent of the choice of the cut
off function $\psi$.
\end{lem}

This is Proposition 5.1 of \cite{dz}.

These two lemmas together  show that \be \label{ko}
\overline{\eta}(X,E, g)= \lim_{a \ra \infty}
\overline{\eta}(D^{\psi,g}_{[0, a]})  \ee for any cut off function
which may depend on $a$ ((\ref{ko}) is to be interpreted as an
equation mod $\mathbb Z$ without the invertibility assumption).
This is exactly the adiabatic limit.
\newline

We now recall the setup and result from \cite{d} on the adiabatic limit of eta invariant, which is an
extension of \cite{bc} to manifolds with boundary. More precisely,
let \be Y\ra M\stackrel{\pi}{\ra} B \ee be a fibration where the
fiber $Y$ is closed but the base $B$ may  have nonempty boundary.
 Let $g_B$ be a metric on $B$ which is of the
product type near the boundary $\pt B$. Now equip $M$ with a
submersion metric $g$,
\[ g = \pi^{*} g_{B} + g_Y  \]
so that $g$ is also product near $\pt M$. This is equivalent to
requiring $g_Y$ to be independent of the normal variable near $\pt
B$, given by the distance to $\pt B$.

The adiabatic metric $g_x$ on $M$ is given by \be g_x =x^{-2}
\pi^{*} g_{B} + g_Y, \ee where $x$ is a positive parameter.

For simplicity we assume that $M$ as well as the vertical tangent
bundle $T^VM$ are spin. Associated to these data we have in
particular the total Dirac operator $D_x$ on $M$, the boundary Dirac
operator $D_x^{\pt M}$ on $\pt M$, and the family of Dirac operators
$D_Y$ along the fibers. If the family $D_Y$ is invertible, then,
according to \cite{bc}, the boundary Dirac operator $D_x^{\pt M}$ is
also invertible for all small $x$, therefore the eta invariant of
$D_x$ with the APS boundary condition,  $\eta(D_x)$, is
well-defined. We have the following result from \cite{d}.

\begin{theo}  \label{al} Consider the fibration $Y\ra M\ra B$ as above. Assume that the
Dirac family along the fiber, $D_Y$, is  invertible. Consider the
total Dirac operator $D_x$ on $X$ with respect to the adiabatic
metric $g_x$ and let $\eta(D_x)$ denote the eta invariant of $D_x$
with the APS boundary condition. Then the limit $\lim_{x
\rightarrow 0} \bar{\eta}(D_x) = \lim_{x \rightarrow 0}
\frac{1}{2} \eta(D_x)$ exists in ${\Bbb R}$ and \be \label{alf}
\lim_{x \rightarrow 0} \bar{\eta}(D_x) =
\int_B\hat{A}\left(\frac{R^B}{2\pi}\right) \wedge \tilde{\eta},
\ee where $R^B$ is the curvature of $g_B$, $\hat{A}$ denotes the
the $\hat{A}$-polynomial  and $\tilde{\eta}$ is the $\eta$-form of
Bismut-Cheeger \cite{bc}.
\end{theo}

Recall that the (unnormalized) $\eta$-form of Bismut-Cheeger, the
$\hat{\eta}$ form, is defined as \be \label{ehf} \hat{\eta} =\left\{
\begin{array}{ll}  {\displaystyle
\int}_{0}^{\infty}{\rm tr}_{s}\left[\left(D_{Y} +
\frac{c(T)}{4t}\right) e^{-B_{t}^{2}}\right]
\frac{dt}{2t^{1/2}} & \mbox{if $\dim Y=2l$}  \\
 & \\
 {\displaystyle \int}_{0}^{\infty} {\rm
tr}^{\rm even}\left[\left(D_Y + \frac{c(T)}{4t}\right)
e^{-B_{t}^{2}}\right] \frac{dt}{2t^{1/2}} & \mbox{if $\dim
Y=2l-1$}
\end{array} \right. ,  \ee
assuming that $\ker D_Y$ does define a vector bundle on $B$. Here
$B_t$ denotes the rescaled Bismut superconnection: \be \label{rbs}
B_t =\tilde{\nabla}^u + t^{1/2}D_Y - \frac{c(T)}{4t^{1/2}}. \ee

We normalize $\hat{\eta}$ by defining \be \label{etf} \tilde{\eta} =
\left\{
\begin{array}{ll} {\displaystyle \sum \frac{1}{(2\pi
i)^j} [\hat{\eta}]_{2j-1}}  & \mbox{if $\dim Y=2l$}  \\
{\displaystyle \sum \frac{1}{(2\pi i)^j} [\hat{\eta}]_{2j}} &
\mbox{if $\dim Y=2l-1$}
\end{array} \right..  \ee

We now turn to the intrinsic spectral interpretation of our eta
invariant.

\begin{theo} \label{mtic} Under the assumption that $\label{ta} \ker [D_X + s\ c(g^{-1} dg)] =0,\ \ \forall\  0
\leq s \leq 1$,
$$
\overline{\eta} \left(X, E,g\right)= \frac{i}{4\pi} \int_0^1
   \int_{0}^{\infty}
{\rm tr}_{s}\left[c\left(g^{-1} dg\right)\left(D_{X} + s\,
c\left(g^{-1} dg\right)\right) e^{-t\left(D_{X} + s\,
c\left(g^{-1} dg\right)\right)^2}
 \right] dt \ ds .
$$
\end{theo}

\Pf We apply Theorem \ref{al} to our current situation where $M=[0,
1] \times X$ fibers over $[0, 1]$ with the fibre $X$. The operator
\[ D_{[0,1]}^{\psi, g} =D^E + (1-\psi)
g^{-1}[D^E, g]=D^E + (1-\psi) c(g^{-1} dg) \] is of Dirac type,
and of product type near the boundaries. Hence the result still
applies.

By the invertibility assumption there is no spectral flow
contribution and hence, by (\ref{ko}), $ \overline{\eta} \left(X,
E,g\right)$ is given by the adiabatic limit formula.

The Dirac family along the fiber is $D_X + (1-\psi(x))c(g^{-1} dg)$.
The curvature of the Bismut superconnection is given by
\[ B_t^2= t \left[D_X + (1-\psi(x))c\left(g^{-1} dg\right)\right]^2 - t^{1/2} \psi'(x) dx\, c\left(g^{-1}
dg\right). \]

Thus,
\[ \hat{\eta} =  \frac{\psi'(x)dx}{2}
\int_{0}^{\infty} {\rm tr}_{s}\left[c\left(g^{-1} dg\right)\left(D_X
+ (1-\psi(x))c\left(g^{-1} dg\right)\right) e^{-t \left(D_X +
(1-\psi(x))c\left(g^{-1} dg\right)\right)^2}\right] dt .\]

Since $\tilde{\eta}= \frac{1}{2\pi i} \hat{\eta}$ here, the
adiabatic limit formula in Theorem \ref{al} gives
\[ \begin{array}{l}
\overline{\eta}(X, E,  g)   \\
  =   \displaystyle\int_0^1
 \frac{\psi'(x)}{4\pi i}
   \int_{0}^{\infty} {\rm tr}_{s}\left[c\left(g^{-1} dg\right)\left(D_X + (1-\psi(x))c\left(g^{-1}
dg\right)\right) e^{-t \left(D_X + (1-\psi(x))c\left(g^{-1} dg\right)\right)^2}\right] dt dx \\
 =  \displaystyle \frac{i}{4\pi } \int_0^1
\int_{0}^{\infty} {\rm tr}_{s}\left[c\left(g^{-1} dg\right)\left
(D_{X} + s\, c\left(g^{-1} dg\right)\right) e^{-t\left(D_{X} + s\,
c\left(g^{-1} dg\right)\right)^2}
 \right] dt \ ds \end{array} \]
 as claimed.
\qed

\section{The noninvertible case}

For a fibration over the circle, Witten's Holonomy Theorem
\cite{w,bf,c} says that the adiabatic limit of the eta invariant
of the total space is related to the holonomy of the determinant
line bundle of the family operators along the fibers. Indeed, in
the invertible case, namely the family operators along the fibers are invertible,
there is an explicit formula for the adiabatic limit of the eta
invariant in terms of the family operators, \cite[(3.166)]{bf},
\cite[(1.56)]{c}, which states \be \label{altic} \lim_{x
\rightarrow 0} \bar{\eta}(D_x) = \frac{i}{4\pi} \int_{S^1}
\int_0^{\infty} {\rm tr}_{s}\left[\left(\tilde{\nabla}^u
D_Y\right)D_Y e^{-tD_Y^2}\right] dt. \ee Of course, the integrand
in the formula (\ref{altic}) is just the degree one term of the
Bismut-Cheeger $\eta$-form.

If one applies (\ref{altic}) to the family $s \in [0,1]
\longrightarrow D_X + s\ c(g^{-1} dg)$, we would obtain Theorem
\ref{mtic}. However, the family here is not periodic. Nevertheless,
it is almost periodic in the sense that the operators at the
endpoints differ by a conjugation. This leads us to the
generalization to the general noninvertible case.

To deal with the noninvertible case, we make use of the framework
and result of \cite{df}. We first recall the setup of \cite{df}.

Suppose $M$ is a compact odd dimensional Riemannian manifold with nonempty
boundary. For simplicity, we assume $M$ is spin so that one can
consider the Dirac operator $D_M$ (the same consideration can be
adapted to Dirac type operators). Further, assume that the metric is of product type near the boundary.
In order to consider eta invariant, one needs to impose boundary conditions and the self adjoint APS boundary condition
amounts to a ``trivialization" of the graded kernel of the boundary Dirac operator $D_{\pt M}$. Taking this into consideration, the result of \cite{df} says that the exponentiated eta invariant of  $D_M$ actually defines an element of the inverse determinant line of the boundary Dirac operator $D_{\pt M}$.

More precisely, let $K_{\pt M}^+\oplus K_{\pt M}^-$ be the (graded) kernel of $D_{\pt M}$ and $\Det^{-1}_{\pt M}$ the inverse determinant line of $D_{\pt M}$:
\be \label{idl}  \Det^{-1}_{\pt M}=\Lambda^{{\rm max}}K_{\pt M}^+\hat{\otimes} [\Lambda^{{\rm max}} K_{\pt M}^-]^{-1} .\ee
Here inverse denotes the dual.

A self adjoint APS boundary condition is determined by a choice of  isometry
 \be \label{sdbc}  T\ :K_{\pt M}^+ \longrightarrow K_{\pt M}^-.  \ee
Let $\reta(D_M(T))$ denote the reduced eta invariant of $D_M$ with
the self adjoint APS boundary condition determined by $T$ (cf.
\cite{df}). A basic result of  \cite{df} says that
\be \label{brodf}  \tau_M = e^{2\pi i \reta(D_M(T))}( \det T)^{-1} \in  \Det^{-1}_{\pt M} \ee
is independent of $T$ (and satisfies the laws of TQFT as well as a variation formula).

Relevant to our discussion here is Witten's Holonomy Theorem as
formulated in this framework. Let $\pi :Y\to Z$ be a fibration whose
typical fiber is a closed even dimensional manifold, and as before
we assume that both $Y$ and $T^VY$ are spin for simplicity. Let
$L\to Z$ denote the corresponding inverse determinant line bundle.
It comes equipped with a (Quillen) metric and a natural unitary
(Bismut-Freed) connection $\nabla $. The curvature of $\nabla $ is
\cite[Theorem 1.21]{bf}
  $$ \Omega ^L = -2 \pi i \left[\int_{Y/Z}\hat{A}(\Omega ^{Y/Z})\right]_{(2)}. $$

Given  $\gamma: [0, 1]\to Z$ a smooth path, let $Y_\gamma =\gamma
^*Y$ denote the pullback of $\pi :Y\to Z$ via $\gamma $; then $\pi
_\gamma :Y_\gamma \to [0, 1]$ is a fibration, the induced fibration.
Let $g_{[0,1]}$ denote an arbitrary metric on the unit interval and
$g_{Y/Z}$ the metric on the vertical tangent bundle $T^VY$. Define a
family of metrics on $Y_\gamma $ by the formula
$$ g_\epsilon = \frac{g_{[0,1]}}{\epsilon ^2}\oplus g_{Y/Z},\qquad \epsilon \not= 0. $$
(We assume that $\gamma$ is constant near the two endpoints so that $g_\epsilon$ is of the product type near the boundary.)

The construction above gives rise to a linear map
\be \label{bcodf}  \tau _{Y_\gamma }(\epsilon ) :L_{\gamma (0)}\longrightarrow L_{\gamma (1)}. \ee

\begin{theo}[Dai-Freed]   The adiabatic limit $\tau _\gamma =\lim\limits_{\epsilon \to0}\tau _{Y_\gamma
}(\epsilon )$ exists and gives the holonomy along $\gamma$ of the Bismut-Freed connection.
\end{theo}

Consider now the fibration $\pi: \mathbb R \times X \longrightarrow \mathbb R$ given by the projection, with the family of Dirac type operators
\be \label{fdo}  s\in \mathbb R \ra D_s = D_X + s c(g^{-1}dg).\ee
Let $L\ra \mathbb R$ be the inverse determinant line bundle with the
Quillen metric and the Bismut-Freed connection. Denote by $L_s$ the fiber of $L$
at $s\in \mathbb R$. Since $D_1=D_X +
c(g^{-1}dg)=g^{-1}D_X g=g^{-1}D_0 g$, there is an isomorphism
\be \label{isom}  g^{-1}: \, L_0 \simeq L_1\ee
determined by the isomorphism $g^{-1}$ between the graded kernels
$\ker D_0$ and $\ker D_1$. On the other hand, since $\mathbb R$ is one dimensional, any
two monotonic paths from $0$ to $1$ are reparametrizations of each
other. Hence there is a unique holonomy map
\be \label{hm}  \tau_{0,1}: L_0 \ra L_1. \ee
Composing with the isomorphism (\ref{isom}) gives rise to a map
\be \label{im}  L_0 \ra L_1 \simeq L_0 \ee
which can then be identified with a complex number $\tau\in \mathbb
C$. In fact, since both the holonomy map (\ref{hm}) and the isomorphism (\ref{isom}) are
unitary maps, $\tau$ has modulus one.

We can now state the main result of this paper as follows.

\begin{theo}  \label{mt} We have
\[ \tau =  e^{2\pi i \reta \left(X, E,g\right)}. \]
\end{theo}

\Pf By taking the exponential we discount the contribution from the
spectral flow to our eta invariant. Thus we are only concerned with
 $\overline{\eta}(D^{\psi, g}_{[0, 1]})$. By definition, $\overline{\eta}(D^{\psi, g}_{[0, 1]})$
is the reduced eta invariant of
\[ D^{\psi, g}=D^E + (1-\psi)
g^{-1}[D^E, g] \] on the cylinder $[0, 1]\times X$  with the
boundary condition $P_{X}(L) $ on one of the boundary components
$\{0\}\times X$ and the boundary condition ${\rm Id}-g^{-1}P_{X}(L)
g $ on the other boundary component $\{1\}\times X$, where $L$ is a
Lagrangian subspace of $\ker D_X$. Let $\ker D_X=K^+_X \oplus K^-_X$
be its $\mathbb Z_2$ grading. Then an isometry
\[ T:\ K^+_X \ra K^-_X \]
gives rise to a Lagrangian subspace, namely the graph of $T$. A
little linear algebra shows that the boundary condition ${\rm
Id}-g^{-1}P_{X}(L) g $ corresponds to the isometry
\[ g^{-1}T^{-1}g:\ g^{-1}K^-_X \ra g^{-1} K^+_X. \]
Hence, we have by the previous theorem and the definition (using the
notation $a$-$\lim$ to denote the adiabatic limit)
\begin{eqnarray*} \tau_{0, 1} & = & a\mbox{-} \lim e^{2\pi i
\overline{\eta}(D^{\psi, g}_{[0, 1]})} (\det T)^{-1} \det g^{-1}T g
 \\
&=& \lim_{a\ra \infty} e^{2\pi i \overline{\eta}(D^{\psi, g}_{[0,
a]})} (\det T)^{-1} \det g^{-1}T g \\
&=& e^{2\pi i \overline{\eta}(D^{\psi, g}_{[0, 1]})} (\det T)^{-1}
\det g^{-1}T g. \end{eqnarray*} Therefore
\[ \tau =  e^{2\pi i \reta \left(X, E,g\right)} \]
using the identification. \qed

$\ $

\begin{rk} Recall that in [DZ, Remarks 2.5 and 5.9], the
$\eta$-invariant is used to give an intrinsic analytic
interpretation of the Wess-Zumino term in the WZW theory. Now by
Theorem 4.2, this term is further interpreted by using holonomy.
\end{rk}

$\ $

\begin{rk} Theorem 4.1 gives an adiabatic limit formula for
(reduced) eta invariants without invertibility assumption for one
dimensional manifolds with boundary, namely the interval. Theorem
3.3, on the other hand, is such a formula with invertibility
assumption, but for any compact manifold with boundary as the
base. It will be interesting to have a general result combining
these two. This will be addressed elsewhere.
\end{rk}

\section{Final remarks}

We end this paper by recalling  a conjecture from \cite{dz}, and
also by   some remarks.

As we mentioned before, the eta type invariant $\overline{\eta}(X,
E, g)$, which we introduced using a cut off function, is in fact
independent of the cut off function. This leads naturally to the
question of whether $\overline{\eta}(X, E, g)$ can actually be
defined directly. The following conjecture is stated in \cite{dz} and \cite{z2}.

Let $D^{[0, 1]}$ be the Dirac operator on $[0, 1]\times X$. We equip
the boundary condition $gP_{X}(L) g^{-1}$ at $\{0\}\times X$ and the
boundary condition ${\rm Id}-P_{X}(L)$ at $\{ 1 \}\times X$.

Then $(D^{[0, 1]}, gP_{X}(L) g^{-1} , {\rm Id}-P_{X}(L)
 )$ forms a self-adjoint elliptic boundary problem. We denote
the corresponding elliptic self-adjoint operator by
$D^{[0,1]}_{gP_{X}(L) g^{-1} , P_{X}(L) }$.

Let $\eta(D^{[0,1]}_{gP_{X}(L) g^{-1}, P_{X}(L)  },s)$ be the
$\eta$-function of $D^{[0,1]}_{gP_{X}(L) g^{-1} , P_{X}(L)
 }$. By [KL, Theorem 3.1],  one knows
that the $\eta$-function $\eta(D^{[0,1]}_{gP_{X}(L) g^{-1} ,
P_{X}(L)  },s)$ admits a meromorphic extension to ${\bf C}$ with
poles of order at most 2. One then defines, as in [KL, Definition
3.2], the $\eta$-invariant of $D^{[0,1]}_{gP_{X}(L) g^{-1}, P_{X}(L)
}$, denoted by $\eta(D^{[0,1]}_{gP_{X}(L) g^{-1} , P_{X}(L)  })$, to
be the constant term in the Laurent expansion of
$\eta(D^{[0,1]}_{gP_{X}(L) g^{-1}, P_{X}(L)  },s)$ at $s=0$.

Let $\overline{\eta}(D^{[0,1]}_{gP_{X}(L) g^{-1}, P_{X}(L)
 })$ be the associated reduced $\eta$-invariant.

$\ $

\noindent {\bf Conjecture}:
$$
\overline{\eta}(X, E, g)=
\overline{\eta}\left(D^{[0,1]}_{gP_{X}(L) g^{-1} , P_{X}(L)
}\right).
$$
\newline

We would also like to say a few words about the technical assumption that
${\rm ind}\, D^E_+=0$ imposed in order to define the eta invariant $\overline{\eta}(X, E, g)$.
The assumption guarantees the existence of the Lagrangian subspaces
$L$ which are used in the boundary conditions. In the Toeplitz index theorem, this
assumption is automatically satisfied since $X=\pt M$ is a boundary. In general, of course,
it may not. However, if one is willing to overlook the integer contribution (as one often does
in applications), this technical issue can be overcome by using another eta invariant, this time
on $S^1 \times X$, as follows. Note that we now have no boundary, hence no need for boundary conditions!

Consider $S^1 \times X= [0, 1] \times X / \sim$ where $\sim$ is the equivalence relation that identifies $0\times X$ with $1\times X$. Let $E_g \ra S^1 \times X$ be the vector bundle which is $E\otimes {\mathbb C}^N$ over $(0, 1) \times X$
and the transition from $0\times X$ to $1\times X$ is given by $g:\ X \ra U(N)$.  Denote by $D_{E_g}$ the Dirac operator  on $S^1 \times X$ twisted by $E_g$.

\begin{prop}\footnote{We thank Jean-Michel Bismut for pointing this out to us several years ago.} One has
\[  \overline{\eta}(X, E, g) \equiv \overline{\eta}(D_{E_g}) \ \ \ {\rm mod} \ \ \mathbb Z . \]
\end{prop}

This is an easy consequence of the so called gluing law for the
eta invariant, see \cite{bu,bl,df}. An analog of this result in
the noncommutative setting plays an important role in \cite{x},
which also contains  an odd dimensional analog of \cite{lmp}.

$\ $

\begin{rk} An application of the Witten holonomy theorem (\cite{w,
bf, c}) to the right hand side of the above formula  leads to an analogous result as
Theorem \ref{mt}. However the family of operators here is not as explicit as in Theorem  \ref{mt}.
\newline
\end{rk}

\begin{rk}  It might be interesting to note the duality that
$\overline{\eta} (X, E,g)$
 is a spectral invariant associated to a $K^1$-representative
on an {\it even} dimensional manifold, while the usual
Atiyah-Patodi-Singer $\eta$-invariant ([APS1]) is a spectral
invariant associated to a $K^0$-representative on an {\it odd}
dimensional manifold.
\end{rk}
$\ $

\end{document}